\documentclass{amsart}

\usepackage{amssymb}

\usepackage{mathrsfs} 

\usepackage[all]{xy} 
\input xy 

\xyoption{all}
\xyoption{2cell} 
\xyoption{v2}

\DeclareMathOperator{\Aut}{Aut}

\DeclareMathOperator{\Out}{Out}

\DeclareMathOperator{\bBun}{\mathbf{Bun}}

\DeclareMathOperator{\id}{id}

\DeclareMathOperator{\Map}{Map}
\DeclareMathOperator{\Ad}{Ad}

\DeclareMathOperator{\Diff}{Diff}

\DeclareMathOperator{\bBibun}{\mathbf{Bibun}}

\DeclareMathOperator{\bBisp}{\mathbf{Bisp}}
\DeclareMathOperator{\bSp}{\mathbf{Sp}}

\DeclareMathOperator{\Type}{Type}

\DeclareMathOperator{\bBGerbe}{\mathbf{BGrb}}

\theoremstyle{plain}
\newtheorem{theorem}{Theorem}[section]
\newtheorem{corollary}[theorem]{Corollary}
\newtheorem{lemma}[theorem]{Lemma}
\newtheorem{proposition}[theorem]{Proposition}

\theoremstyle{definition}
\newtheorem{definition}[theorem]{Definition}

\theoremstyle{remark}
\newtheorem{example}{Example}[section]
\newtheorem{note}{Note}[section]

\numberwithin{equation}{section}
\numberwithin{figure}{section}

\newcommand{\cG}{\mathcal{G}}

\newcommand{\ZZ}{{\mathbb Z}}

\renewcommand{\a}{\alpha}

\begin{document}

\title[Bibundles]{On the existence of bibundles}
  \author[M. Murray]{Michael Murray}
  \address[Michael Murray]
  {School of Mathematical Sciences\\
  University of Adelaide\\
  Adelaide, SA 5005 \\
  Australia}
  \email{michael.murray@adelaide.edu.au}

 \author[D.M. Roberts]{David Michael Roberts}
  \address[David Michael Roberts]
  {School of Mathematical Sciences\\
  University of Adelaide\\
  Adelaide, SA 5005 \\
  Australia}
  \email{david.roberts@adelaide.edu.au}
  
\author[D. Stevenson]{Danny Stevenson}
\address[Danny Stevenson]
{School of Mathematics and Statistics\\
University of Glasgow  \\
15 University Gardens\\
Glasgow G12 8QW\\
United Kingdom}
\email{Danny.Stevenson@glasgow.ac.uk} 
   
  \thanks{The first author acknowledges the support of the Australian
Research Council. The second author 
acknowledges the support of an Australian Postgraduate Research Award.  The third author was supported by 
the Engineering and Physical Sciences Research Council [grant number EP/I010610/1].  
Finnur Larusson and Klaas Landsman
are thanked for useful comments made after talks on this topic. We thank the referee for many useful comments.}



\begin{abstract}
{We consider the existence of bibundles, in other words locally trivial principal $G$ spaces with commuting
left and right $G$ actions. We show that their existence is closely related to the structure
of the group $\Out(G)$ of outer automorphisms of $G$.  We also develop a classifying theory for bibundles.  The theory is developed
in full  generality for $(H, G)$ bibundles for a crossed-module $(H, G)$ and we show with 
examples the close links with loop group bundles. }
\end{abstract}
\maketitle

\tableofcontents

\section{Introduction}
There has been  interest recently in gerbes that have non-abelian band \cite{AshCanJur}  particularly 
for applications to string theory.   In the literature that has developed it is notable that  it seems to be difficult 
to find concrete examples which are 
not closely related to abelian gerbes.   A closer examination suggests that the problem centres around
the need, when defining a gerbe, to be able to form a product of principal $G$ bundles in such a way as to yield another 
principal $G$ bundle (rather than a $G \times G$ bundle).  A similar problem arises in module theory when wants to take 
a tensor product of modules of a non-commutative ring $R$.  In this case one is led to consider bimodules, i.e.\ modules  over the ground ring with commuting left and right actions of $R$.    By analogy, we are naturally led to study not just principal bundles but so-called `bibundles'.  
These are  fibrings that 
are simultaneously left and right principal $G$ bundles in such a way that the left and right $G$ actions commute. 

However, the existence of such objects is more problematic.  To see why, 
consider the fibre of such a bibundle. In the case of a $G$ bundle the fibre is a right $G$ space and there is only one of these up to isomorphism. In the case of bibundles the 
fibre is a $G$ {\em bispace} and now there are different isomorphism classes labelled by $\Out(G)$, the group of 
outer automorphisms of $G$. We dwell in some apparently pedantic detail on the structure of $G$ bispaces.  This effort 
however is rewarded  by making many constructions for bibundles  immediate. 

Fundamental to our approach is the idea of the {\em type} of a bispace or bibundle.  In the case of 
bispaces the type of a bispace is an element of $\Out(G)$ which classifies its isomorphism class. In the case of 
a  bibundle $P \to M$ it is a map from $M $ to $\Out(G)$ whose value at $m \in M$ classifies the isomorphism
class of the bispace which is the fibre of $P$ over $m$.   We call the map which associates to a bibundle the type 
of each of its fibres the {\em type map}. It forms part of an exact sequence of pointed sets
\begin{equation*}
 \pi_0\bBun_{Z(G)}(M) \stackrel{\iota}\longrightarrow \pi_0\bBibun_G(M) \stackrel{\Type}{\longrightarrow} \Map(M, \Out(G))) 
\end{equation*}
where $Z(G)$ is the center of $G$, and $\pi_0\bBun_{Z(G)}(M)$ is the pointed set of isomorphism classes of $Z(G)$-bundles over $M$.
The message that follows from the exactness of this sequence of pointed sets is that for genuinely non-abelian bibundles to exist we need   $\Out(G)$ to be large. 
In the case that $G$ is simple and simply-connected $\Out(G)$ is well known to be 
the (small) finite group of automorphisms of the Dynkin diagram of $G$.  More interesting examples arise when $G$ is the group $\Omega K$
of based loops in a compact group $K$ whose outer automorphism group has large subgroups such as $K$ 
itself. 

In summary then we start in Section 2 with a detailed discussion of  bispaces. As we will see, bispaces are 
a partial `categorification' of the notion of $G$-spaces in which the structure group is replaced by a certain kind of groupoid --- 
a so-called `2-group' or crossed module.  For simplicity 
in this introduction we have only considered $G$ bispaces, which correspond to a restricted class of 2-groups.  
To obtain a more flexible theory we will need to discuss the more
general case of $(H,G)$ bispaces for a crossed module $(H, G)$. In  Section 3 we consider the case
of $(H, G)$ bibundles and as well as explaining the exact sequence above we consider the classifying theory 
of bibundles. Finally in the conclusion we indicate briefly some results on the more complicated case of
$(H, G)$ bibundle gerbes.

\section{Bispaces}
\label{sec:bispaces}

Let $G$ be a topological group and let $X$ be a space.  If $G$ acts freely and transitively on both the left and the 
right of $X$ and  these actions commute, we  call  $X$ a  $G$ {\em bispace}. We will very often be interested in a smooth version of this notion 
where $G$ is a Lie group, $X$ is a manifold and both of the actions of $G$ on $X$ 
are smooth --- for convenience we will also call these objects $G$ bispaces.  
 However, it turns out that it is more natural to consider a more general notion, that of an $(H, G)$ bispace, where $(H, G)$ 
 is a `crossed module'.  In the next subsection we will give some motivation for this, building up to the definition of $(H,G)$ bispace (see Definition~\ref{def:(H,G) bispace} below).  

\subsection{$(H,G)$ bispaces}

First consider a $G$ bispace $X$ and,  following \cite{AshCanJur,Breen1}, define a map
$$
\psi \colon X \to \Aut(G)
$$
by $xg = \psi(x)(g)x$. Note that $\psi(x)$ is indeed in $\Aut(G)$ because we have
$\psi(x)(gh) x = x gh =\psi(x)(g) x h = \psi(x)(g)\psi(x)(h) x$. 
We call $\psi$ the {\em structure map} of $X$. We have 
$$
\psi(xg)(h) x g  = x gh = x (ghg^{-1}) g = \psi(x)(ghg^{-1}) xg = \psi(x) \Ad(g) (h) xg
$$
so that $\psi(xg) = \psi(x) \circ \Ad(g)$.   Thus $\psi$ is right $G$ equivariant if we consider $X$ to be a right 
$G$ space and $\Aut(G)$ a right $G$ space under the adjoint action.  The data of the right $G$ space $X$ 
together with the equivariant map $\psi\colon X\to \Aut(G)$ is sufficient to recover the bispace $X$, more 
precisely we have the following lemma from \cite{Breen1}.

\begin{lemma}[\cite{Breen1} Lemme 2.5] 
\label{lemma:structure}
The structure map of a bispace gives rise to  an equivalence between
\begin{enumerate}
\item $G$ bispaces $X$ 
\item Pairs $(X, \psi)$ consisting of a right $G$ space $X$ and an 
equivariant map $\psi \colon X \to \Aut(G)$.
\end{enumerate}
\end{lemma}

A slightly more general idea would be to choose a subgroup $H  \subset \Aut(G)$ containing $\Ad(G)$ and require that the structure map $\psi$ 
take values in $H$, in other words $\psi(x)\in H$ for all $x\in X$. 
We will take one step beyond this.  Recall (see for instance \cite{BaezLauda, Breen1}) that a 
 {\em crossed module} is a pair of topological groups $(H, G)$ together with homomorphisms 
$$
G  \stackrel{t}{\to} H  \stackrel{\alpha}{\to}  \Aut(G)
$$
satisfying the following two conditions:
\begin{enumerate}
\item $t$ is $H$-equivariant   for the action of $H$ on $G$ defined by $\alpha$ and the adjoint action $\Ad_H$ of $H$ on $H$, 
that is $t(\a(h)(g)) = h t(g) h^{-1}$, 
\item the action of $G$ on itself induced by $t$ is the adjoint action of $G$ on itself, i.e.\  $\alpha \circ t = \Ad_G$.
\end{enumerate}

Note that we have the following example.

\begin{example}
If $H \subset\Aut(G)$ is a subgroup containing $\Ad(G)$ then 
$$
G  \stackrel{\Ad}{\to} H \hookrightarrow \Aut(G).
$$
is a crossed-module. 
\end{example}

A {\em morphism} of crossed modules $(H,G)\to (H',G')$ consists of a pair of homomorphisms $u\colon H\to H'$ and $v\colon G\to G'$ such that 
the diagram 
$$ 
\xymatrix{ 
G \ar[d]_t \ar[r]^-{v} & G' \ar[d]^-{t'} \\ 
H \ar[r]_-{u} & H' } 
$$ 
commutes and the equivariance condition $v(\a(h)(g)) = \a'(u(h))(v(g))$ is satisfied.  

Two easy consequences of the definition of crossed module are the following: $G_1 = \ker(t) \subset Z(G)$, the centre of $G$ (hence $\ker(t)$ is abelian), and also $t(G) \subset H$ is normal.  
Therefore we have exact sequences of groups 
\[ 
1\to G_1\to G\to G/G_1\to 1
\] 
and
\[ 
1\to t(G)\to H\to H/t(G)\to 1. 
\]
Throughout  we will assume that the projections $G\to G/G_1$ and $H\to H/t(G)$ admit local sections (if $H\to G$ is a crossed module of Lie groups this is of course immediate, unless one is 
dealing with some classes of infinite dimensional Lie groups).  
We will 
sometimes adopt the notation $G \to H$ for a crossed module.  The following definition appears in \cite{Breen1}.

\begin{definition}[\cite{Breen1} page 432]
\label{def:(H,G) bispace}
Let $(H, G)$ be a crossed module.  An $(H, G)$ {\em bispace} is a pair $(X, \psi)$ consisting of a right  $G$ space $X$ and an equivariant map $\psi \colon X \to H$. 
\end{definition}

By equivariant we mean that $\psi(xg) = \psi(x) t(g)$. We shall often write bispace instead of $(H, G)$ bispace
when the context is clear and  call the map $\psi$ the {\em structure map} of $X$. Note that our definition of $(H, G)$ bispace is different to that in \cite{AshCanJur} where such a thing, in our notation, is a left $H$-space with an equivariant map to $\Aut(G)$.

\begin{example}
If $G$ is a topological group then there is a canonical crossed module
$$
G  \stackrel{\Ad}{\to} \Aut(G) \stackrel{\id}{\to}  \Aut(G).
$$
It is straightforward from Lemma \ref{lemma:structure} to see that an $(\Aut(G), G)$ bispace is the same thing as a $G$ bispace.
\end{example}

If $X$ is an $(H, G)$ bispace with structure map $\psi \colon X \to H$
then $\a \circ \psi \colon X \to \Aut(G)$ is equivariant and we have the following lemma.
\begin{lemma}
\label{lem:(H,G)bispace = Gbispace}
If $(X, \psi)$ is an $(H, G)$ bispace then $(X, \a \circ \psi)$ is a $G$-bispace (via the correspondence 
of Lemma~\ref{lemma:structure}).
\end{lemma}
This means that an $(H, G)$ bispace has a left action of $G$ defined by 
\begin{equation}
\label{eq:left}
g x = x ((\a\circ \psi) (x))^{-1}(g)
\end{equation}
which we use in the future, often without comment.  Note however that $\psi(gx) = t(g)\psi(x)$.  

 If $h \in H$ denote by $[h]$ the coset in the quotient group $H/t(G)$.  By equivariance $\psi$ defines a unique element 
$\phi = [\psi(x)] \in H/t(G) $ which we call the {\em  type} of $X$.   If $1$ denotes a point we have the commuting diagram
\begin{equation*}
\label{diagram}
\begin{array}{ccc}
X  & \stackrel{\psi}{\longrightarrow}& H\\
   \downarrow &   &  \downarrow \\
1 & \stackrel{\phi}{\longrightarrow} & H/t(G)
   \end{array}
\end{equation*}
where $\phi(1) = \phi$.  We write $\Type(X)$ for the type of $X$.

In the case of a $G$ bispace $X$ we have that $\Type(X) \in \Aut(G)/\Ad(G) = \Out(G)$, the group 
of outer automorphims of $G$. We have seen (Lemma~\ref{lem:(H,G)bispace = Gbispace}) 
that any $(H, G)$ bispace is also a $G$ bispace
and hence has a type in $\Out(G)$. This is the image of the type in $H/t(G)$ under the 
homomorphism
$$
H/t(G) \to \Out(G)
$$
induced by $\a$.

\begin{example}
\label{ex:trivial-twists}
Choose an element $\xi \in H$ and let $X = G$, considered as a right $G$-space under
group multiplication. Define a structure map $\psi\colon G \to H$  by $\psi(x) = \xi t(x)$. 
The induced bispace left action on $X$ is given by 
$$
k \star x  =  x (\a \psi(x))^{-1}(k)
          = x (\a( \xi t(x) ))^{-1}(k) 
           = x (\a(\xi) \Ad(x) )^{-1}(k).
$$
Denote this $(H, G)$ bispace by $T(\xi)$. Then
$\Type(T(\xi)) = [\xi]$. 

In particular we have the {\em trivial} bispace $T = T(1)$ whose structure map is
$t \colon G \to H$ and for which the induced $G$ bispace is just $G$ with the usual left and 
right $G$ action.
\end{example}

\begin{example}
\label{ex:abelian-bispace}
Let $X$ be a right $A$ space for an abelian group $A$. We can make $X$ an $A$ bispace by defining $a x b = xba^{-1}$ and
with this definition $X$ has structure map $\psi$ defined by 
$\psi(x)(a) = a^{-1}$.  
Note that 
$$
A \to 1 \to \Aut(A)
$$
is a crossed module precisely when $A$ is abelian. A right $A$ space is then a $(1, A)$ bispace with structure map 
equal to $1$.  It is also possible to consider $(\Aut(A), A)$ bispaces such as  the Jandl bundle
discussed below in Example \ref{ex:jandl}.
\end{example}

If $X$ and $Y$  are $G$ bispaces then a function $f \colon X \to Y$ is called a $G$ bispace {\em morphism} if it commutes with the left and right actions.  Note that, just as in the case of ordinary $G$-spaces, every morphism 
of $G$-bispaces is an isomorphism.  We have the following lemma (see Remarque 2.7 of 
\cite{Breen1}).    

\begin{lemma}
If $X$ and $Y$ are $G$ bispaces and $f \colon X \to Y$ is a  bijection then $f$ is a $G$ bispace isomorphism
if and only if it commutes with the right action and the structure map.
\end{lemma}
\begin{proof}
Follows from Lemma \ref{lemma:structure}.
\end{proof}

It is then natural to make the following definition.

\begin{definition}
If $X$ and $Y$  are $(H, G) $ bispaces then a function $f \colon X \to Y$ is called an $(H, G)$ bispace {\em morphism} if it commutes with the right action and the structure maps.  
\end{definition}

Denote by $\bBisp_{(H,G)}$ the category of all $(H,G)$ bispaces and bispace morphisms.    In the special case of the crossed module $(\Aut(G),G)$ associated to a 
group $G$, we will denote the corresponding category of bispaces by $\bBisp_G$.  Note that 
all $(H,G)$ bispace morphisms are automatically 
isomorphisms and hence    
$\bBisp_{(H,G)}$ is in fact a groupoid.  
We have the following proposition.  

\begin{proposition}
\label{prop:type inj}
Two $(H, G)$ bispaces $X$ and $Y$ are isomorphic if and only if they have the same type. 
\end{proposition}
\begin{proof}
Denote by $\psi \colon X \to H$ and $\chi \colon Y \to H$  the structure maps of $X$ and $Y$ respectively.  If $f \colon X \to Y$ is an isomorphism 
then clearly $\chi = \psi \circ f$ and hence the type of $Y$ is equal to the type of $X$. 

Conversely assume the types  are both equal to $\xi \in H/t(G)$. Notice that $\chi$ and $\psi$ are onto the preimage of $\xi$ in $H$  so we can choose $x\in X$ and $y\in Y$ such that $\psi(x) = \chi(y)$. 
 Define  $f \colon X \to Y$ by $f(xg)  = yg$.  Then $f$ is a 
bijection and commutes with the right $G$ action by construction.  Moreover the $G$ equivariance of the structure maps shows
that $\chi \circ f = \psi$ giving the required result.
\end{proof}

We can interpret this result as saying that there is a functor $\Type\colon \bBisp_{(H,G)}\to H/t(G)$, where $H/t(G)$ is considered as a {\em discrete} groupoid, i.e.\ there are no non-identity morphisms.  

The groupoid $\bBisp_{(H,G)}$ has extra structure: there is a functor 
$$
\otimes\colon \bBisp_{(H,G)}\times \bBisp_{(H,G)}\to \bBisp_{(H,G)}
$$
 which sends a pair of bispaces $(X,Y)$ to the {\em product} bispace 
$X\otimes Y$ which is defined as follows.  
If $(X,\psi_X)$ and $(Y,\psi_Y)$ are $(H, G)$ bispaces 
then $X\otimes Y$ is defined to be the bispace $X\otimes Y = (X\times Y)/G$, where $G$ acts on $X\times Y$ by  
\[
(x, y)g = (xg, g^{-1} y) = (xg,  y(\alpha\psi_Y(y))^{-1}(g^{-1})).
\] 
Denote the equivalence class of $(x, y)$ in $X\otimes Y$ by $x \otimes y$. There is a natural 
  right action of $G$ on $X\otimes Y$ given by $(x\otimes y)g = x\otimes (yg)$.  
  Define a map $\psi \colon X \otimes Y \to H$  by $\psi(x \otimes y) = \psi_X(x) \psi_Y(y)$.   
It is straightforward to check that this is well-defined and a 
structure map for $X \otimes Y$ making $(X\otimes Y,\psi)$ 
an $(H, G)$ space. The left action induced by this structure map can be calculated using equation 
\eqref{eq:left} to be $g(x \otimes y) =  gx \otimes y$. 

 It is also straightforward to check that the 
process of forming products of bispaces is functorial.  
Just as is the case when forming tensor products of modules, 
the product of bispaces is not strictly associative, 
however it is associative up to a canonical natural isomorphism.  
Note also that the type of the product bispace satisfies
$\Type(X \otimes Y) = \Type(X) \Type(Y)$.

\begin{example}
\label{ex:twisting}
If $\xi \in H$ we  denote $X \otimes T(\xi)$ by $X(\xi)$.  It is straightforward to show that 
$[x, g] \mapsto x \a(\xi)(g)$ defines a bijection from $X \otimes T(\xi)$ to $X$ with 
inverse $x \mapsto [x, 1]$. With this identification the type map is $x \mapsto [\psi(x) \xi]$
and the right action is $x g = x \a(\xi)(g) $. 
\end{example}

There is also a functor $(-)^*\colon \bBisp_{(H,G)}\to \bBisp_{(H,G)}$ which sends a bispace $X$ to its {\em dual} $X^*$.  
If $(X,\psi)$ is an $(H, G)$ bispace 
we define the dual $(H, G)$ bispace $(X^*,\psi^*)$ to be the same 
underlying space  $X$, but with the structure map $\psi^* = \psi^{-1}$ and the right group action 
$x \cdot g = x \a(\psi^{-1}(x))(g^{-1})$. Again, using equation \eqref{eq:left} it can be seen that $g \cdot x = g^{-1}  x$. 

The following lemma is straightforward.
\begin{lemma}
\label{lemma:canonical-iso}
For any $(H, G)$ bispace $X$ we have canonical isomorphisms 
$X\otimes T \cong T\otimes X\cong X$ and $X\otimes X^*$ isomorphic to $T$.  
\end{lemma}

Write $\pi_0\bBisp_{(H,G)}$ for the set of isomorphism classes in the groupoid $\bBisp_{(H,G)}$.  This is a pointed set, pointed by the isomorphism class of the trivial bispace.  The functors $\otimes$ and 
$(-)^*$ induce a corresponding product and notion of dual on $\pi_0\bBisp_{(H,G)}$.  Since the process of forming products of bispaces is associative up to a canonical natural isomorphism the product on 
$\pi_0\bBisp_{(H,G)}$ is associative.  Lemma~\ref{lemma:canonical-iso} shows that in fact $\pi_0\bBisp_{(H,G)}$ has the structure of a group.  

The fact that the type map for bispaces preserves products means that the functor $\Type\colon \bBisp_{(H,G)}\to H/t(G)$ induces a homomorphism of groups 
\[ 
\Type\colon \pi_0\bBisp_{(H,G)}\to H/t(G).  
\]
Proposition~\ref{prop:type inj} shows that this map is injective, and Example~\ref{ex:trivial-twists} shows that it is surjective.  Hence $\Type$ is an isomorphism of groups.  We summarize this discussion in the 
next proposition.  
\begin{proposition}
The type map induces an isomorphism of groups 
\[
\Type\colon \pi_0\bBisp_{(H,G)} \cong H/t(G). 
\]
\end{proposition}

\subsection{Extension and reduction of bispaces}
\label{sec:ext-red}

Recall that we denote the kernel of $t \colon G \to H$ by $G_1$.  We have seen above that $G_1\subset Z(G)$.  Notice that if $h \in H$ and $g \in G_1$ then  
$t ( \a (h)(g )) = h t(g) h^{-1} = 1$ so that $\a(h)(G_1) \subset G_1$. Hence there is an action 
of $H$ on $t(G) = G/G_1$ via $\a$ and a crossed module
$$
t(G)  \hookrightarrow H  \stackrel{\alpha}{\to}  \Aut(t(G)).
$$
We have the following lemma. 
\begin{lemma}
If $X$ is an $(H, G)$ bispace then $X/G_1$, with the right $G$-action 
and structure map induced from $X$, is an $(H, t(G))$ bispace.
\end{lemma}

If $G \to H$ is a crossed module we say a  {\em crossed submodule} is a crossed module $G_0 \to H_0$ with the property that 
$G_0$ is a subgroup of $G$, $H_0$ is a subgroup of $H$, $t(G_0) \subset H_0$ and the elements in $\a(H_0) \subset \Aut(G)$
fix $G_0$ and thus define a homomorphism $H_0 \to \Aut(G_0)$. In such a case $G_0 \to H_0 $ is clearly a 
crossed module and the inclusions define a  morphism of crossed modules. 

\begin{example} 
If $t \colon G \to H$ is a crossed module and $G_1 = \ker(t)$ then $(1, G_1)$ is a crossed submodule
of $(H, G)$ where $1$ is the identity subgroup of $H$.
\end{example}   

Let $X$ be an $(H, G)$ bispace and $(H_0, G_0)$ a crossed-submodule of $(H, G)$. We say that $X_0 \subset X$
is a {\em reduction} of $X$ to $(H_0, G_0)$ if $X_0$ is an orbit of $G_0$ and $\psi(X_0)\subset H_0$. 
Clearly $X_0$ is an $(H_0, G_0)$ bispace with structure map $\psi_{X_0}$. 

Let $X_1$ denote the subspace
\begin{equation}
\label{eq:reduction}
X_1 = \{ x \in X \mid \psi(x) = 1 \}.
\end{equation}
Then we have the following lemma. 

\begin{lemma}
\label{lemma:reduction}
Let $X$ be an $(H, G)$ bispace then 
\begin{enumerate}
\item $X_1$ is non-empty if and only if $\Type(X) = 1$. 
\item If $X_1$ is non-empty then it is a reduction of $X$ to $(1, G_1)$.
\end{enumerate}
\end{lemma}
\begin{proof}
The first part is obvious from the definition. 
If $X_1$ is  non-empty and  if $ x \in X_1$ and 
$g \in G_1$ then $\psi(x g) = \psi(x) t(g) = 1$.  On the other hand if $x, y \in X_1$ then $x = yg$
for some $g \in G$ and $1 = \psi(x) = \psi(y) t(g) = t(g)$ so that $g \in G_1$. 
\end{proof}

\begin{example}
If $X$ is a $G$ bispace then $t = \Ad \colon G \to  \Aut(G) = H$ so that $G_1 = Z(G)$. 
Then $X_1$ is either empty or a  $(1,Z(G))$ bispace. 
\end{example}

Let  $(H, G)$ and $(H', G')$ be crossed modules,  $(\zeta, \eta) \colon (H, G) \to (H', G')$ a morphism
of crossed modules and $X$ an $(H, G)$ bispace. Let $X(G') = X \times_G G'$ where  the action of 
$G$ on $X \times G'$ is $(x, g') g = (xg, \eta(g^{-1}) g')$. Clearly $X(G') $ is a right $G'$
space with action $[x, g']g'' = [x, g' g'']$. Define a structure map by $\psi([x, g']) = 
\zeta(\psi(x))t'(g')$. This is well-defined and equivariant making $X(G')$ an $(H', G')$ bispace. 
We call $X(G')$ the {\em extension} of $X$ to $(H', G')$. 
In particular we have the following Lemma whose proof is straightforward.

\begin{lemma} 
\label{lemma:extension}
Let $X_0$ be a reduction of the $(H, G)$ bispace $X$ to the crossed submodule $(H_0, G_0)$. 
Then the map $ X_0(G) \to X$ defined by $[x, g] \mapsto xg$ defines an isomorphism of
$(H, G)$ bispaces.
\end{lemma}

Thus we have also:
\begin{lemma}
If $\Type(X) = 1$ then $X \simeq X_1(G)$.
\end{lemma}

\begin{note}  As additional motivation for our introduction of crossed modules we note that if we have a 
$G$ bispace $X$ and a homomorphism $G \to H$ there is no natural induced $H$ bispace.  We need the additional data
of a homomorphism $\Aut(G) \to \Aut(H)$ so that the two homomorphisms give rise to a homomomorphism of 
crossed modules $(\Aut(G), G) \to (\Aut(H), H)$.  Indeed this is the key reason for considering crossed modules as
coefficient objects for nonabelian gerbes - we don't have
functoriality in $G$ when considering the corresponding cohomology
theory, only functoriality for maps of crossed modules.
\end{note}

\subsection{2-groups and crossed modules}

\label{2gp}
Bispaces and crossed modules are closely related to 2-groups.  Recall (see for instance \cite{BaezLauda, Breen1}) that a (topological) {\em 2-group} is a groupoid object $\mathcal{G}$ in the category of topological groups.  
We will not spell out what this means precisely: suffice it to say that it means that $\mathcal{G}$ is a groupoid for which both the objects and morphisms have the structure of topological groups.  
Since we will only ever be interested in topological 2-groups we will omit the adjective `topological'.  

Every crossed module $t\colon G\to H$ uniquely determines a 2-group $\cG$ and vice versa.  The group of objects of $\cG$ is defined to be $H$, while the group of 
morphisms of $\cG$ is defined to be the semi-direct product $H\rtimes G$.  For more details we refer to \cite{BaezLauda, Breen1}.  

The notion of 2-group that we have surveyed here is a `strict' one; there is also a notion of {\em weak} 2-group for which we refer to \cite{BaezLauda}.  Put briefly, a weak 2-group is a (topological) 2-groupoid with one object; from 
another point of view a weak 2-group consists of a groupoid $\mathcal{G}$ equipped with a functor $\otimes\colon \cG\times \cG\to\cG$ which is associative up to a coherent natural isomorphism, together with several other axioms.  

The groupoid $\bBisp_{(H, G)}$ of all $(H, G)$ bispaces is a prime example of a weak 2-group in this sense.  As we have seen, the product $X\otimes Y$ of two $(H, G)$ bispaces defines a functor 
$$
\otimes \colon \bBisp_{(H, G)}\times \bBisp_{(H, G)} \to \bBisp_{(H, G)}
$$
 as above.  There is a canonical functor from the groupoid $\bBisp_{(H,G)}$ to the groupoid $\cG$ associated to the crossed module $(H,G)$.  This canonical functor preserves products in $\bBisp_{(H,G)}$ and 
 $\cG$ in an appropriately weakened sense: it turns out that this canonical functor is an {\em equivalence} between the weak 2-group $\bBisp_{(H, G)}$ and its strict version $\cG$.  For more details we refer to \cite{BaezLauda}.  
 
 This last statement is partially analogous to the following well known fact about topological groups: if $G$ is a topological group and $\bSp_G$ denotes the groupoid of right $G$ spaces and maps between them, then 
 $\bSp_G$ is equivalent to $G$, thought of as a groupoid with one object.  The difference between this situation and the one we are considering lies in the fact that the groupoid $G$ is not normally a 2-group, in fact it is a 2-group if and only if 
 $G$ is abelian.  
 
 The groups $\ker(t)$ and $H/t(G)$ also have a nice interpretation in terms of 2-groups.  It turns out (using the technology of simplicial homotopy theory, see for instance \cite{MaySOAT}) that one can make sense of the homotopy groups 
 $\pi_i(\cG)$ of a 2-group $\cG$.  In fact the crossed module $G\to H$ associated to a 2-group $\cG$ arises in this setting as the {\em Moore complex} of the simplicial group which is the nerve of $\cG$.  If one follows the standard recipe for 
 computing the simplicial homotopy groups of a simplicial group then one finds that the homotopy groups $\pi_i(\cG)$ correspond to the homology groups of $G\to H$, thought of as a complex concentrated in degrees 0 and 1.  So one finds that 
 $\pi_0(\cG) = H/t(G)$ and that $\pi_1(\cG) = \ker(t)$.  The fact that $t(G)$ is normal in $H$ and the fact that $\ker(t)$ is abelian can then be understood as higher dimensional analogues of the fact that the set of path components of a topological group has the structure of a group and the 
 fact that the fundamental group of a topological group is abelian, respectively.       

\section{Bibundles}

\begin{definition}
Let $(H, G)$ be a crossed module. If  $P \to M$ is a  (right) principal $G$  bundle with an equivariant map 
$\psi \colon P \to H$ such that each fibre of $P \to M$ is a $(H, G)$ bispace we 
call $P \to M$ an $(H, G)$ {\em bibundle}.  
\end{definition}

We will call $P$ the {\em total space} and $M$ the {\em base space} of an $(H,G)$ bibundle $P\to M$.  
As for the case of bispaces we say that a {\em morphism} of bibundles is a morphism of the underlying principal bundles 
which commutes with the structure maps.  Clearly every morphism of bibundles inducing 
the identity on base spaces is an isomorphism.  We will write 
$\bBibun_{(H,G)}(M)$ for the groupoid of bibundles on $M$ and we will denote the set of isomorphism classes of 
bibundles on $M$ by $\pi_0\bBibun_{(H,G)}(M)$.  If $P$ is a bibundle on $M$ then we will write $[P]$ for its isomorphism 
class in $\pi_0\bBibun_{(H, G)}(M)$.  

Consider an $(H, G)$ bibundle $P \to M$. Each fibre of $P \to M$ is an $(H, G)$ bispace so  it follows immediately from the discussions in Section \ref{sec:bispaces}  that we have a commuting diagram
\begin{equation}
\label{eq:diagram}
\begin{array}{ccc}
 P  & \stackrel{\psi}{\longrightarrow}& H\\
   \downarrow &   &  \downarrow \\
M & \stackrel{\phi}{\longrightarrow} & H/t(G) 
   \end{array}
\end{equation} 
where $\psi$ satisfies $pg = (\alpha\circ\psi)(p)(g)p$.  As before we will call $\psi$ the {\em structure map} of $P$ and 
$\phi$ the {\em type map} of $P$.  Local triviality of $P \to M$ will ensure that $\psi$ and $\phi$ are smooth 
or continuous as appropriate.

\begin{note}
We remark that the notion of $(H,G)$ bibundle (and the notion of $(H,G)$ bispace) is actually a special case of the notion of groupoid bundle.  Recall that if $\cG$ is a topological groupoid 
with space of objects $G_0$ and space of morphisms $G_1$, then a $\cG$ {\em groupoid bundle} over $M$ (see \cite{HS,Moer}) consists of a map 
$\pi\colon P\to M$ admitting local sections together with an {\em action} of $\cG$ on $P$, in other words the data of 
\begin{enumerate} 
\item a map $p:P\to G_0$, 

\item a map $m\colon P\times_{G_0} G_1\to P$ 
\end{enumerate} 
satisfying certain axioms (for which refer to \cite{HS, Moer}).  Furthermore the action is required to be principal in the sense that 
the diagram 
$$ 
\xymatrix{ 
P\times_{G_0} G_1 \ar[d]_-{p_1} \ar[r]^-m & P \ar[d] \\ 
P \ar[r] & M } 
$$ 
is a pullback, where $p_1$ denotes projection onto the first factor.  When $\cG$ is the groupoid associated to a crossed module $(H,G)$ as described in Subsection~\ref{2gp} above, the notion of $\cG$-bundle coincides exactly with the notion of 
$(H,G)$-bibundle.    
\end{note}
 
\begin{example}
\label{ex: abelian bibundle}
If $P \to M$ is a (right) principal $A$ bundle for an abelian group $A$ then $P \to M$ is a  $(1, A)$ bibundle  with 
structure map $1 \colon P \to 1$. 
\end{example}

\begin{example}
\label{ex:groups}
Suppose that $G$ is a  normal subgroup of $H$ so that we have a crossed module 
$i\colon G\to H$.  Let $K$ denote the quotient group $H/G$ and suppose that the projection $H \to K$ admits local sections,  then 
  $H \to K$ is an $(H, G)$ bibundle.  In this case the structure map $\psi\colon H\to H$ is the identity.  
    \end{example}

  \begin{example}
Similarly if $(H, G)$ is  a crossed module then $H \to H/t(G)$ is an $(H, t(G))$ bibundle with structure map $\psi \colon H \to H$ equal to the identity.  Notice  that   the fibre $H_\xi$ over $\xi \in H/t(G)$  is a $t(G)$ bispace of type $\xi$. 
\end{example}

  \begin{example}
  \label{ex:loopgroup}
  As an important example of the above construction let $K$ be a simple, simply-connected, compact Lie group and denote by $PK$ the group of all smooth maps $k \colon [0, 1]  \to K$ with $k(0) = 1$. If we 
  define $\pi \colon PK \to K$ to be evaluation of a path at $1$ then this is an $\Omega K$ bibundle. Here  we are defining the loop group $\Omega K \subset PK $ to be  the subgroup of all paths $k$ with $k(0) = k(1)$. 
  Note that this is a larger group than the group of smooth maps $k$ from $S^1$ to $K$ with  $k(0) = 1$. As in the general case above  the adjoint action of $PK$ on itself fixes the subgroup $\Omega K$ so we have
  a crossed module
  $$
  \Omega K \to P K \to \Aut(\Omega K) 
  $$
  and thus the $(PK, \Omega K)$ bibundle $PK \to K$.
    
  \end{example}

\begin{example}
If $P \to M$ is an $(H, G)$ bibundle and 
$G_1$ is the kernel of $t \colon G \to  H$ then  $P/G_1 \to M$ is a $(H,t(G))$ bibundle, where $t(G) = G/G_1$. 
\end{example}

Just as with the structure and type maps many of  the other notions we have introduced for bispaces can be extended immediately to bibundles
by applying them to the fibres of $P \to M$. In particular this applies to the notions of reduction and extension and the product and dual constructions. So if $P \to M$ and $Q \to M$ 
are $(H, G)$ bibundles then there are bibundles $P^* \to M$ and $P \otimes Q \to M$. 
If $P$, $Q$ and $R$ are bibundles on $M$ then there are canonical isomorphisms $P\otimes (Q\otimes R)\cong (P\otimes Q)\otimes R$ and $P\otimes P^*\cong T$ where $T = M\times G$ is the 
trivial bibundle on $M$ whose fibre at each point of $M$ is the trivial bispace.  In a completely analogous way to the earlier discussion for bispaces, we have the following Lemma.  

\begin{lemma} 
\label{group structure on IBibun} 
The set $\pi_0\bBibun_{(H, G)}(M)$ of isomorphism classes of bibundles on $M$ forms a group with product $[P]\otimes [Q]$ defined by $[P\otimes Q]$ and where the inverse $[P]^{-1}$ of an 
element $[P]$ is given by $[P^*]$.   
\end{lemma} 

\begin{example} 
\label{ex: group bibundle morphism}
If $H\to K$ is the bibundle of Example~\ref{ex:groups} above then the product $H\times H\to H$ in $H$ induces a bibundle morphism $H\otimes H\to H$ covering the product 
in $K$ so that the diagram 
$$ 
\xymatrix{ 
H\otimes H\ar[d] \ar[r] & H \ar[d] \\ 
K\times K \ar[r] & K }
$$ 
commutes.  Similarly the inverse map $h\mapsto h^{-1}$ in $H$ defines an isomorphism $H\to H^*$ covering the inverse map in $K$.    
\end{example}

Note that \eqref{eq:diagram} is  not quite a morphism of $(H, G)$ bibundles, instead 
\begin{equation*}
\begin{array}{ccc}
 P/G_1  & \stackrel{\psi}{\longrightarrow}& H\\
   \downarrow &   &  \downarrow \\
M & \stackrel{\phi}{\longrightarrow} & {H}/t(G)
   \end{array}
\end{equation*}
is a morphism of $(t(G), H)$ bibundles.  
From this discussion we deduce the following proposition.

\begin{proposition} If $P \to M$ is an $(H, G)$ bibundle with type map 
$\phi \colon M \to H/t(G)$ then $P/G_1$ is the pull-back of $H \to H/t(G)$ by $\phi$.
\end{proposition}

As a consequence we can deduce the following corollary.

\begin{corollary} If $t \colon G \to H$ is a crossed module with $\ker(t) = 1$ then $G \simeq t(G)$ and 
every bibundle $P \to M$ is the pullback of the $(H, G)$ bibundle $H \to H/t(G)$
by the type map.
\end{corollary}

Consider the case when the type map is the constant map to the identity in $H/t(G)$.  Then 
each fibre of $P \to M$ has a non-empty subset of points $p\in P$ such 
that $\psi(p) =1 $. By analogy with \eqref{eq:reduction} above Lemma \ref{lemma:reduction} 
denote  the union of these subsets by $P_1$ and note
that from  Lemma \ref{lemma:reduction} we have that $P_1$ is  a reduction of $P$ to $(1, t(G))$. 
 
 Following  \cite{AshCanJur} we say that a section $s$ of $P$ is a {\em central}
 section if $\psi \circ s = 1$.  We then have the following proposition.

\begin{proposition}[c.f \cite{AshCanJur}]
\label{prop:trivial1}
A bibundle $P$ is trivial if and only if it has a central section. 
\end{proposition}

This gives us immediately the following result.

\begin{corollary}
\label{cor:trivial2}
A bibundle $P$ is trivial if and only if the type map is equal to $1$ and $P_1$  is trivial. 
\end{corollary}

 If $Q \to M$ is a $G_1$ bundle, that is a $(1, G_1)$ bibundle,  then we can extend to an  $(H, G)$ bibundle 
 $\iota(Q) = Q(G)$ using the construction from Subsection \ref{sec:ext-red}.
 So we have  $\iota(Q) = Q \times_{G_1} G$ and  $\psi([q, g]) = t(g) \in H$.
 The type of $\iota(Q)$ is clearly $1 \in \Map(M, H/t(G))$. On the other hand assume that $P$ has 
 type $1$ so that we have a well-defined reduction of $P_1$ to $(1,G_1)$.
 Then from Lemma \ref{lemma:extension} we have the isomorphism 
$$
\begin{array}{ccc}
 P_1 \times_{G_1} G & \simeq & P \\{}
[p, g] & \mapsto & pg .
\end{array}
$$
Hence we have the following proposition.
\begin{proposition}
$Q \to M$ is a trivial $G_1$ bundle if and only if $\iota(Q)$ is a  trivial $(H, G)$ bibundle.
\end{proposition}
\begin{proof}
Clearly if $Q$ has a section then it induces a section of $\iota(Q)_1$ which
is a central section on $\iota(Q)$ so by Corollary  \ref{cor:trivial2} $Q$ is 
trivial. On the other hand if $\iota(Q)$ is trivial then  it has a 
central section. But that must be a section of $\iota(Q)_1 \simeq Q$, thus $Q$ has a section 
and is trivial.
\end{proof}

It follows that we have a sequence of 2-groups and homomorphisms between them
\begin{equation}
\label{2gp exact seq}
1\to \bBun_{G_1}(M)\stackrel{\iota}\longrightarrow \bBibun_{H,G}(M) \xrightarrow{\Type} \Map(M, H/t(G)),
\end{equation}
where the group $\Map(M, H/t(G))$ is thought of as a discrete 2-group.  This sequence is `exact' in the 
following sense.  The homomorphism $\iota$ is faithful, and if $P$ is a bibundle on $M$ then 
$\Type(P) = 1$ if and only if $P$ is isomorphic to 
a bibundle of the form $\iota(R)$, where $R$ is a $G_1$ bundle on $M$.  

On passing to isomorphism classes 
we obtain the exact sequence of sets 
\begin{equation}
\label{eq:bibundle_seq}
1\to \pi_0\bBun_{G_1}(M)\to \pi_0\bBibun_{H,G}(M)\xrightarrow{\Type} \Map(M,H/t(G)).  
\end{equation}
Consider now the image of the type map. To understand this let  $\phi \colon M \to H/t(G)$ be any map.  We can 
pullback the $t(G)$ bundle $H \to H/t(G)$ along this map so that we get a pullback diagram
\begin{equation*}
\begin{array}{ccc}
 \phi^*(H)  & \stackrel{\psi}{\longrightarrow}& H\\
   \downarrow &   &  \downarrow \\
M & \stackrel{\phi}{\longrightarrow} & {H}/{t(G)}.
   \end{array}
\end{equation*}
Assume we can find a $G$ bundle $P \to M$ which lifts the (right) $t(G)$  bundle $ \phi^*(H) \to M$
to a (right) $G$ bundle. Then we can define a  map $\psi \colon P \to H$
by the composite $P \to  \phi^*(H) \to H$. This map $\psi$ is equivariant and hence 
is the structure map for an $(H, G)$ bibundle structure on $P \to M$ with $\phi$ as type map.  We conclude
that for $\phi$ to be in the image of $\Type$ it suffices for us to be able to 
lift $\phi^*(H)$ to a $G$ bundle.  Consider then the central extension 
$$
0 \to G_1 \to G \to t(G) \to 0.
$$
The obstruction to lifting $\phi^*(H)$ from $t(G)$ to $G$ is the non-triviality of the 
class in $H^2(M,G_1)$ of the $G_1$ lifting bundle gerbe  \cite{Murray} associated to $\phi^*(H)$.
It follows that we have the exact sequence of groups 
\begin{equation}
\label{exact pi0 seq}
1\to \pi_0\bBun_{G_1}(M) \stackrel{\iota}{\longrightarrow} 
\pi_0\bBibun_{H,G}(M) \xrightarrow{\Type} \Map(M, H/t(G)) \stackrel{\epsilon}\longrightarrow H^2(M,G_1).  
\end{equation}
We remark that there is an alternative way to arrive at this exact sequence, for which we sketch the details.  
The 2-category $\bBGerbe_{G_1}(M)$ of $G_1$ bundle 
gerbes on $M$ together with the corresponding 1-morphisms and 2-morphisms between them, is an example of a 
3-group.  The map which sends a map $\phi\colon M\to H/t(G)$ to the lifting $G_1$ bundle gerbe on $M$ 
determined by the pullback $G_1$ bundle $\phi^*H$ defines a homomorphism $\Map(M,H/t(G))\to \bBGerbe_{G_1}(M)$.
We can extend the exact 
sequence~\eqref{2gp exact seq} one term to the right and obtain an exact sequence of 3-groups 
\[
1\to \bBun_{G_1}(M)\stackrel{\iota}\longrightarrow \bBibun_{H,G}(M) \xrightarrow{\Type} \Map(M, H/t(G))
\to \bBGerbe_{G_1}(M) 
\]
where `exact' is to be understood in a similar sense to that above.  
Taking $\pi_0$ recovers the exact sequence~\eqref{exact pi0 seq} above.  

The exact sequence \eqref{eq:bibundle_seq} tells us loosely that if $H/t(G)$ is small then most $(G, H)$ bibundles
are likely to be abelian, i.e.\  reduce to abelian $(1, G_1)$ bibundles or $G_1$ bundles.  Recall that there is a map $H/t(G) \to \Aut(G)/\Ad(G) = 
\Out(G)$ so the question of whether or not there are many $(G, H)$ bibundles that are not abelian also relates to the size of $\Out(G)$
which we consider next.

\subsection{Type maps that lift}
We say that a map $\phi \colon M \to H/t(G)$ {\em lifts} if there is some $\hat\phi \colon M \to H$
which projects to $\phi$.  The construction used in Example \ref{ex:trivial-twists} can be applied to define
a bibundle $T(\hat\phi)$ whose fibre over $x \in M$ is $T(\hat\phi(x))$. This gives a map 
$$
\Map(M, H) \to \pi_0\bBibun_{(H, G)}(M)
$$
which makes the following diagram commute
\begin{equation} 
\label{lifting} 
\xy 
(-20,7.5)*+{\Map(M, H)}="1"; 
(20,7.5)*+{\pi_0\bBibun_{(H, G)}(M)}="2"; 
(20,-7.5)*+{\Map(M, H/t(G).}="3"; 
{\ar "1";"2"}; 
{\ar "1";"3"}; 
{\ar^{\Type} "2";"3"}
\endxy
\end{equation}
Recall also from Example \ref{ex:twisting} that if $Q \to M$ is a  bibundle we denote 
$Q(\hat\phi) = Q \otimes T(\hat \phi)$.

As a result we have the following proposition.
\begin{proposition} 
If $P \to M$ is a bibundle with type map $\phi \colon M \to H/t(G)$ which lifts to $\hat\phi \colon M \to H$
 then $P$ is isomorphic  to $\iota(R)(\hat\phi)$ for some $G_1$ bundle $R$.
\end{proposition}
\begin{proof} 
By construction $T(\hat\phi)$ has type map $\phi$, the same as $P$. Hence there is a $G_1$  bundle $R$ with 
$P \simeq \iota(R)(\hat\phi)$.
\end{proof}

In particular  we deduce the following when $G$ is a simply-connected and semi-simple Lie group.
\begin{proposition} 
\label{prop:lift}
If  $G$ is a simply-connected and semi-simple Lie group then every $G$ bibundle is of the form 
$\iota(R)(\hat\phi)$ for some $\hat\phi \colon M \to \Aut(G)$.
\end{proposition}
\begin{proof}
We have that $H/t(G) = \Aut(G)/\Ad(G) = \Out(G)$ is the group of automorphisms of the Dynkin diagram of $G$ and
hence discrete. It follows that  the type map $\phi$ is constant on connected components of $M$
and hence  lifts.
\end{proof}

\begin{example}
\label{ex:jandl}
Consider an $(\Aut(A), A)$ bibundle where $A$ is abelian.  In this case $\Ad(A) = 1$ so that $\phi$ lifts and the lift is, in fact, just $\widehat \phi = \phi$.  
Hence from Proposition \ref{prop:lift} every bibundle has the form $P(\widehat\phi)$ for some $A$-bundle $P$.   A particular case is $A = U(1)$ when $\Aut(U(1)) = \ZZ_2$. In this case we have a $U(1)$ bundle $P \to M$ and on each connected component we give it a left action by defining $z p = p z^{\pm 1}$ depending on whether $\phi$
restricted to that connected component is $\pm 1$.  A $(\ZZ_2, U(1))$ bibundle 
is a {\em Jandl bundle} \cite{NikSch}.
\end{example}

\begin{note}
Just as principal $G$-bundles on $M$ can be described in terms of local data via their transition cocycles $g_{ij}\colon U_i\cap U_j\to G$ relative to some open cover $\{U_i\}$ of $M$, so also do bibundles have such a local description.  
If $P$ is an $(H,G)$ bibundle on $M$ for some crossed module $(H,G)$ then, for a sufficiently fine open cover $\{U_i\}$ of $M$, we can associate to $P$ families of maps $g_{ij}\colon U_i\cap U_j\to G$ and $h_i\colon U_i\to H$ satisfying the 
cocycle conditions 
\begin{gather*}
g_{ij}  g_{jk}  = g_{ik} \\
h_j  = h_i t(g_{ij}). 
\end{gather*}
The maps $g_{ij}$ are the usual transition cocycles of the bundle $P$ and arise from comparing trivializations of $P$ on overlapping patches.  The maps $h_i$ are formed by composing the section defining a trivialization of $P$ over $U_i$ with the  structure map $P\to H$. 
One can introduce an equivalence relation on pairs $(g_{ij},h_i)$ and form a cohomology group $H^0(M,\cG)$ which parametrizes isomorphism classes of $(H, G)$ bibundles.  

\end{note}

\subsection{The bundle of bibundle structures on a principal bundle}
\label{sec:bibundle_structures}

Note that the structure map  $\psi$ of an $(H, G)$ bibundle can be viewed as a section of the associated bundle $P \times_G H$ 
where  $G$ acts on the right of both $P$ and $H$, the latter through the action $h\cdot g = ht(g)$.

Another way of thinking about this is that if $P$ is a right $G$ bundle then each  element in $P \times_G H$ over $x \in M$ defines a bispace structure on the right $G$ space $P_x$. This is because each such 
element defines a structure map $P_x \to H$: 
if $[p, h] \in P \times_G H$ then define a map $P_x \to H$ which
sends $pg \mapsto h t(g)$. Clearly this is well defined and by construction it is 
equivariant.  
 
Note that bibundles pullback along maps: if $P$ is an $(H,G)$ bibundle on $M$ and $f\colon N\to M$ 
is a map, then $f^*P = N\times_M P$ has a natural structure of an $(H,G)$ bibundle on $N$.   
Therefore, if we let $\pi \colon P \times_G H \to M$
be the projection map, then $\pi^*(P) \to P \times_G H$ has a natural structure of an $(H, G)$
bibundle. We can identify $P \times H$ with $\pi^*(P) \subset P \times P \times_G H $ by the map $(p, h)  \mapsto (p, [p, h])$ 
and hence induce an $(H, G)$ bibundle structure on $P \times H \to P \times_G H $. 
This is given by 
$$
(p, h) g = (pg , ht(g))
$$
and the structure map is the projection to $H$.  Thus we have a commutative diagram

\begin{equation*} 
\xy 
(-15,7.5)*+{P \times H}="1"; 
(10,7.5)*+{H}="2"; 
(-15,-7.5)*+{P\times_G H}="3"; 
(10,-7.5)*+{H/t(G).}="4"; 
{\ar^{\Psi} "1";"2"}; 
{\ar "1";"3"}; 
{\ar_-{\Phi} "3";"4"}; 
{\ar "2";"4"};
\endxy
\end{equation*}

In summary then we see that a $(H, G)$ bibundle structure on $P$ is a section $M \to  P\times_G H$ and the pullback via
this section of $P \times H \to P \times_G H$ is naturally isomorphic to $P \to M$ as an $(H, G)$ bibundle.

\subsection{The universal case}

Let $G$ be a topological group and let $EG\to BG$ be the universal $G$ bundle. 
We can apply the construction of Section \ref{sec:bibundle_structures} to the right $G$ space $EG$ and form the space
$$ 
 EG \times H.
$$ 
As we saw this is a bibundle over $EG \times_G H$ with  right $G$ action given by 
\begin{equation}
\label{eq:rightaction}
(e, h)g  = (eg, h t(g))
\end{equation}
and structure map the projection onto $H$.  It follows that the type map $\phi \colon EG \times_G H \to H/t(G)$
sends $[e, h]$ to the equivalence class of $h$ in $H/t(G)$ and we have the commuting
diagram 

\begin{equation}
\label{universal commuting diagram} 
\xy 
(-15,7.5)*+{EG \times H}="1"; 
(10,7.5)*+{H}="2"; 
(-15,-7.5)*+{EG \times_G H}="3"; 
(10,-7.5)*+{H/t(G).}="4"; 
{\ar^{\Psi} "1";"2"}; 
{\ar "1";"3"}; 
{\ar_-{\Phi} "3";"4"}; 
{\ar "2";"4"};
\endxy
\end{equation}

If $P \to M$ is an $(H, G)$  bibundle then there is a classifying map $f \colon M \to BG$ which lifts 
to a right $G$ equivariant map $\hat f \colon P \to EG$. Together with the structure
map $\psi \colon P \to H$ this defines a  homomorphism of $(H, G)$ bibundles 
$$
\begin{array}{ccccc}
\hat \Gamma &\colon &P & \to  & EG \times H \\
& & p &\mapsto & (\hat f(p), \psi(p)). 
 \end{array}
 $$
Denote by $\Gamma \colon M \to EG \times_G H$ the induced map. Then we have a commuting diagram
of bibundles
\begin{equation} 
\label{comm square} 
\xy 
(-15,7.5)*+{P}="1"; 
(10,7.5)*+{EG \times H}="2"; 
(35,7.5)*+{H}="3"; 
(-15,-7.5)*+{M}="4"; 
(10,-7.5)*+{EG \times_G H}="5"; 
(35,-7.5)*+{H/t(G).}="6"; 
{\ar^-{\hat \Gamma} "1";"2"}; 
{\ar^{\Psi} "2";"3"}; 
{\ar_-{\Gamma} "4";"5"}; 
{\ar_-{\Phi} "5";"6"};
{\ar "1";"4"}; 
{\ar "2";"5"}; 
{\ar "3";"6"}; 
\endxy
\end{equation}
As a concrete example of this we have the following example.
\begin{example}
Consider the crossed module $\Omega K \to PK$ so that $G = \Omega K$ and 
$H = PK$.  We can realise $EG \to BG$ as the path fibration $PK \to K$. 
Then $EG \times H = PK \times PK $ with the right action of $\Omega K$. The diagram above becomes

\begin{equation} 
\label{loop comm square} 
\xy 
(-15,7.5)*+{P}="1"; 
(10,7.5)*+{PK\times PK}="2"; 
(35,7.5)*+{PK}="3"; 
(-15,-7.5)*+{M}="4"; 
(10,-7.5)*+{PK \times_{\Omega K} PK}="5"; 
(35,-7.5)*+{K.}="6"; 
{\ar^-{\hat \Gamma} "1";"2"}; 
{\ar^-{\Psi} "2";"3"}; 
{\ar_-{\Gamma} "4";"5"}; 
{\ar_-{\Phi} "5";"6"};
{\ar "1";"4"}; 
{\ar "2";"5"}; 
{\ar "3";"6"}; 
\endxy
\end{equation}
\end{example}
The function  $\Gamma$ in~\eqref{comm square} is a candidate for a classifying map for bibundles. However, before we can justify this we need to 
study how bibundles induced by pullback are related to one another under homotopies of maps.

\subsection{Homotopy of bibundles} 
\label{subsec:homotopy}
A key fact in the theory of ordinary bundles is the fact that 
homotopic maps induce isomorphic bundles under pullback.  
This same statement fails to be true for bibundles.  To see 
why, recall that the type 
map $M\to H/t(G)$ of an $(H,G)$ bibundle $P$ on 
$M$ is an invariant of the isomorphism class of $P$.  In other 
words, if $P$ and $Q$ are isomorphic bibundles on $M$ then the type map 
of $P$ is {\em equal} to the type map of $Q$.  
If $h\colon M\times I\to N$ is a homotopy between maps 
$f_0,f_1\colon M\to N$ and $P$ is a bibundle on $N$ then 
there is no reason why the type map of the induced bundle 
$h^*P$ should be 
constant in the $t$-direction (here $t$ is the coordinate on the interval $I$), and hence there is no reason 
why the bibundles $f_0^*P$ and $f_1^*P$ should be isomorphic.  

To get around this problem we clearly need to restrict 
our attention to homotopies $h\colon M\times I\to N$ which satisfy the following 
property: the composite map $\phi \circ h\colon M\times I\to H/t(G)$ 
is constant in the $t$-direction, where $\phi$ is the 
type map of $P$.  Put another way, 
we have a commutative diagram 
\[
\xy
(-10,5)*+{M\times I}="1";
(10,5)*+{N}="2";
(0,-5)*+{H/t(G)}="3";
{\ar^-h "1";"2"};
{\ar "1";"3"};
{\ar^-{\phi} "2";"3"};
\endxy
\]   
where the map $M\times I\to H/t(G)$ is the composition $M\times I\stackrel{\mathrm{pr_1}}{\to} 
M\stackrel{\phi_0}{\to} H/t(G)$, where 
$\phi_0$ denotes the type map of $h_0^*P$.  

Therefore, we should regard $h\colon M\times I\to N$ as a map in the category of spaces over $H/t(G)$, 
which places us in the realm of parametrized homotopy theory.
  
Recall that if $B$ is a space then the category of spaces over $B$ is the category whose objects 
are spaces $X$ equipped with a map $X\to B$ (we will refer to such spaces as spaces {\em over} $B$) and 
whose morphisms are maps of spaces which are compatible with the projections to $B$ (we will 
also refer to such maps as maps {\em over} $B$).  Note that the identity map $1_B\colon B\to B$ 
exhibits $B$ as an object of this category (this is the terminal object of the category 
of spaces over $B$).  Note also that product  in the category of spaces over $B$ 
of two spaces $X$, $Y$ over $B$ is the 
fibre product $X\times_B Y$.

A homotopy between two maps $f_0,f_1\colon M\to N$ in the 
category of spaces over $B$ is then the usual 
sort of map $h\colon M\times I\to N$, but 
where $h_t \colon M \to N$ defined  by $h_t(x) = h(x, t)$ is a map over $B$ for all $t \in I$.
We will sometimes say that $h$ is a homotopy from $f_0$ to $f_1$ {\em over} $B$.  
Clearly if $f_0$ and $f_1$ are homotopic maps over $B$, then $f_0$ and $f_1$ are homotopic in the usual sense.  

In particular we have the notion of a homotopy equivalence over $B$, also called a {\em fibre 
homotopy equivalence}.  Notice that if 
$f\colon X\to Y$ is a fibre homotopy equivalence, then for any space $Z$ over $B$, the induced 
map $Z\times_B X\to Z\times_B Y$ is also a fibre homotopy equivalence (this fails to be true if $f$ is just an 
ordinary homotopy equivalence).  In particular if $X$ is contractible as a space over $B$, in other words 
if the given map $X\to B$ is a homotopy equivalence over $B$, then for any space $Z$ over $B$ the 
induced map $Z\times_B X\to Z$ is a fibre homotopy equivalence.       

Homotopy of maps over $B$ is an equivalence relation and 
we will write $[X,Y]_{B}$ for the set of homotopy classes of maps over 
$B$ between two spaces $X$ and $Y$ over $B$.  

With this understanding of the notion of homotopy, we have the following proposition.  
\begin{proposition}
\label{prop:homotopy}
Suppose that $f_0\colon M\to N$ and $f_1\colon M\to N$ are homotopic maps in the category of spaces over $H/t(G)$, where $M$ is paracompact.  If $P$ is a bibundle on $N$ then the 
bibundles $P_0 = f_0^*P$ and $P_1 = f_1^*P$ on $M$ are isomorphic.  
\end{proposition}

\begin{proof}
Let $h\colon M\times I\to N$ be a homotopy from $f_0$ to $f_1$ over $H/t(G)$ and let $Q\to M\times I$ denote the induced bundle $h^*P$.  
 It suffices to show that the bibundles $Q$ and $P_0\times I$ on $M\times I$ are isomorphic.  Standard bundle theory shows that the underlying  
 $G$ bundles $P$ and $P_0\times I$ are isomorphic.  

The $G$ bundle $P_0\times I$ is equipped with the structure map 
$P_0\times I\to P_0\to H$ which is constant in the $I$-direction (the map $P_0\to H$ is the restriction 
of the structure map $Q\to H$ of $Q$ to $Q|_{M\times \{0\}}$).  Using the isomorphism $Q\cong P_0\times I$ 
we can define on $P_0\times I$ a new structure map $P_0\times I\to H$.  
We would like to show that these two structure maps 
define isomorphic bibundle structures on $P_0\times I$.  

Let $\hat{\phi}_1\colon P_0\times I\to H$ denote 
the structure map which is constant in the $I$-direction and let $\hat{\phi}_2\colon P_0\times I\to H$ denote the structure
 map which is induced by the isomorphism of $G$-bundles $Q\cong P_0\times I$.  Since $\hat{\phi}_1$ and 
$\hat{\phi}_2$ correspond to the same (constant) type map $M\times I\to H/t(G)$ we must have that 
$\hat{\phi}_2 =\hat{\phi}_1\chi$ for some map $\chi\colon P_0\times I\to G/\ker(t)$ which satisfies 
$$ 
R_g^*\chi = t(g)^{-1}\chi t(g).
$$ 
Suppose that we can find a map $\hat{\chi}\colon P_0\times I\to G$ satisfying $\chi = t(\hat{\chi})$ and 
$R_g^*\hat{\chi} = g^{-1}\hat{\chi} g$.  Then $\hat{\chi}$ defines a $G$ bundle automorphism of $P_0\times I$.  
Furthermore this bundle automorphism is compatible with $\hat{\phi}_1$ and $\hat{\phi}_2$ in the obvious sense. 
It follows that $\hat{\phi}_1$ and 
$\hat{\phi}_2$ define isomorphic bibundle structures on $P_0\times I$.  

It remains to prove the existence of the equivariant map $\hat{\chi}\colon P_0\times I\to G$.  We are given the 
equivariant map $\chi\colon P_0\times I\to G/\ker(t)$ and we note that this restricts to the constant map 
$1$ on $P_0\times \{0\}$.  Since $\chi$ is an equivariant map 
it can be regarded as a section $s$ of the associated bundle of groups $P_0(G/\ker(t))\times I$  
on $M\times I$.   The condition that $\chi$ restricts to the constant map on $P_0\times \{0\}$ translates into the condition 
that the section $s$ is identically $1$ on $M\times \{0\}$.   

The homomorphism $t\colon G\to G/\ker(t)$ induces a map 
$P_0(G)\to P_0(G/\ker(t))$ of the associated bundles.  As part of our assumptions on the crossed 
module $(H,G)$ (see Section~\ref{sec:bispaces}) we assume that the map 
$G\to G/\ker(t)$ has local sections, and is hence locally trivial with fibre $\ker(t)$.  It follows that 
$P_0(G)\to P_0(G/\ker(t))$ is a locally trivial fibre bundle with fibre $P_0(\ker(t))$.     

We need to prove that the section $s$ of $P_0(G/\ker(t))\times I$ lifts to 
a section $\hat{s}$ of $P_0(G)\times I$.  Such a section $\hat{s}$ can be thought of as a section of the bundle $R$ 
on $M\times I$ obtained by pulling back the bundle $P_0(G)\times I\to P_0(G/\ker(t))\times I$ with the section 
$s\colon M\times I\to P_0(G/\ker(t))\times I$.  

We are in the following situation: we have a 
fibre bundle $R\to M\times I$ together with a section defined over 
$M\times \{0\}$ and we want to extend this section to a section 
defined over the whole of $M\times I$.  Since $R\to M\times I$ is a fibration, it has the homotopy 
lifting property, hence such a section exists.  
\end{proof}

\subsection{Classifying theory for bibundles} 

We can consider the total space $P$ of an $(H,G)$ bibundle on $M$ as an object in the category of spaces over $H$ via the structure map.  
Let us say that $P$ is {\em contractible} as a space over $H$ if the structure map $P\to H$ is a homotopy equivalence in the category of spaces over $H$.  As an example, since $EG$ is contractible as an ordinary space, 
$EG\times H$ is contractible when viewed as a space over $H$ with the map to $H$ being 
projection to the second factor.  We have the following analogue of the classical bundle classification theorem.  
\begin{theorem}
\label{classif thm}
Suppose that $E$ is an $(H,G)$ bibundle over a space $B$ such that $E$ is contractible when viewed as a space over $H$.  Then $E\to B$ is a universal bibundle in the sense that 
there is an isomorphism 
\[ 
[M,B]_{H/t(G)}\cong \pi_0\bBibun_{(H,G)}(M)
\]
induced by pullback of bibundles, for any paracompact space $M$.  
\end{theorem}
Our proof will be an adaptation of the proof of Theorem 7.5 from \cite{Dold}.  In this paper Dold 
introduces some key notions which we will recall here, as they 
play an important role in what follows.  
Recall (see Definition 2.2 of \cite{Dold}) that a map 
$p\colon Y\to X$ is said to have the {\em section extension property} if the following is true: 
if $A$ is a closed subspace of $X$ and $s\colon A\to Y$ 
is a section of $p$ defined over $A$ which has an 
extension to a `halo' around $A$, then there is an extension 
of $s$ to a section of $p$ defined on $X$.  Here a 
{\em halo} of $A$ is a subset $V$ of $X$ such that 
there is a continuous map $\tau\colon X\to I$ with the property that 
$\tau(a) = 1$ for all $a\in A$ and $\tau(x) = 0$ for all 
$x\in X-V$.  It follows that any map $p\colon Y\to X$ 
with the section extension property admits at least one section.    

Dold proves the important Theorem 2.7 of \cite{Dold} 
which says that if $\{U_i\}$ is a numerable open cover 
of $X$ and the restriction $p|_{U_i}$ of $p$ to 
$U_i$ has the section extension property for all $i\in I$, 
then $p\colon Y\to X$ has the section extension property.  
A sufficient condition for a map $q\colon U\to V$ to have 
the section extension property is that $q$ is {\em shrinkable}, 
in other words $q$ is fibre homotopy equivalent to the identity 
map $1_V$ (see Proposition 2.3 of \cite{Dold}).  
Therefore, if there exists a numerable open cover $\{U_i\}$ of 
$X$ such that $p|_{U_i}$ is shrinkable for all $i\in I$, then 
$p\colon Y\to X$ has the section extension property.        

\begin{proof}
Proposition~\ref{prop:homotopy} shows that the map 
$[M,B]_{H/t(G)}\to \pi_0\bBibun_{(H,G)}(M)$ which sends a 
homotopy class $[f]$ of maps over $H/t(G)$ to the isomorphism 
class of the pullback bibundle $f^*E$ 
is well defined.  We want to show that this map is an isomorphism.  
We first show that this map is surjective.  Let $P$ be a bibundle on 
$M$ and consider the space $P\times_H E$ over $M$ and its 
quotient $(P\times_H E)/G$ by the diagonal $G$-action.  
A section of the canonical 
map $(P\times_H E)/G\to M$ is a $G$-equivariant map $P\to E$ 
which is compatible with the structure maps of $P$ and $E$.  
If such a section exists, then we have a map $f\colon M\to B$ and it follows that $P$ is isomorphic to  
the pullback $f^*E$.  

First suppose that $(P,\phi)$ is a bibundle on $M$ such that the underlying principal $G$ bundle 
is trivial --- suppose that $s\colon M\to P$ is a section of the map $\pi\colon P\to M$.  We will 
show that in this case the map $(P\times_H E)/G\to M$ is shrinkable.  We regard $M$ as a 
space over $H$ via the map $\phi\circ s\colon M\to H$.  Then there is an isomorphism 
\begin{align*} 
(P\times_H E)/G\to M\times_H E \\ 
[p,e]\mapsto (\pi(p),e\tau(p,s\pi(p))) 
\end{align*}
which is compatible with the projections to $M$.  Here $\tau\colon P^{[2]}\to G$ denotes 
the usual map which satisfies $p_2 = p_1\tau(p_1,p_2)$ for $(p_1,p_2)\in P^{[2]}$.  It is 
easy to see that this map is well defined.     
To see that $(\pi(p),e\tau(p,s\pi(p)))\in M\times_H E$ as claimed, note that if 
$(p,e)\in P\times_H E$ is a representative of $[p,e]$ then $\phi(p) = \psi(e)$, where 
$\psi\colon E\to H$ is the structure map of the bibundle $E$.  Therefore 
\[ 
\psi(e\tau(p,s(m))) = \phi(p)t(\tau(p,s(m))) = \phi(p\tau(p,s(m))) = \phi(s(m)),  
\]
as required.  An inverse for the map above is given by the map $M\times_H E\to 
(P\times_H E)/G$ which sends $(m,e)\mapsto [s(m),e]$ as is easily checked.    
Since $E$ is contractible as a space over $H$, we see that $M\times_H E\to M$, and hence 
$((M\times G)\times_H E)/G\to M$, is shrinkable.  

Since $M$ is paracompact, we can form a numerable cover $\{U_i\}_{i\in I}$ of $M$ by open sets 
$U_i$ with the property that there exist local sections $s_i\colon U_i\to P$ 
of $\pi\colon P\to M$ over $U_i$.  It follows from the previous argument 
that the restriction of $(P\times_H E)/G\to M$ to $U_i$ is shrinkable for all $i\in I$.  Hence, 
by Theorem 2.7 of \cite{Dold}, we see that $(P\times_H E)/G\to M$ admits a section over $M$.  
Hence $P$ is induced by pullback from $E$ via a map $M\to B$.  

Next we show that the map is injective.  Suppose that $f_0\colon M\to B$ and 
$f_1\colon M\to B$ are representatives of homotopy classes of maps 
for which there is a bibundle isomorphism $f_0^*E\cong f_1^*E$.  Let $P_0 = f_0^*E$, $P_1 = f_1^*E$ and suppose that $\a\colon P_0\to P_1$ is a bibundle isomorphism.  Let $P = P_0\times I$ and consider, as in \cite{Dold}, 
the bibundle map 
\begin{equation} 
\label{partial section}
P|_{M\times \{0,1\}}\to E 
\end{equation} 
induced by $\a$ and the bibundle map $\a_0\colon P_0\to E$.  This bibundle map is a section of $(P\times_H E)/G$ over $M\times \{0,1\}$.  By the same argument 
as above, since $E$ is contractible as a space over $H$, the map $(P\times_H E)/G\to M\times I$ is locally shrinkable and hence has the section extension property.  

The section of $(P\times_H E)/G\to M\times I$ over $M\times \{0,1\}$ defined by~\eqref{partial section} has an extension to a halo around $M\times \{0,1\}$ since it can be 
extended (following \cite{Dold}) to a bundle map 
$$
\begin{array}{ccc}
 P|_{M\times ([0,\frac{1}{2})\cup (\frac{1}{2},1])}  &\to& E \\ 
 (u,t)  &\mapsto &\begin{cases} 
\a_0(u)\ & \text{if}\ t<\frac{1}{2}\\ 
\a_1\a(u)\ & \text{if}\ t>\frac{1}{2}
\end{cases} 
\end{array}
$$
where we have written $\a_1\colon P_1\to E$ for the bundle map covering map $f_1\colon M\times I\to B$.  The set $V = M\times ([0,\frac{1}{2})\cup (\frac{1}{2},1])$ is a halo around $M\times \{0,1\}$ as explained in 
\cite{Dold}.  Therefore, by the section extension property, there is an extension of the section~\eqref{partial section} to a global section defined over the whole space $M\times I$.  
This global section corresponds to a bundle map $P_0\times I\to E$ which covers a map $M\times I\to B$.  The latter map is a homotopy between $f_0$ and $f_1$.  
\end{proof}
\begin{note}
Note that this theorem need not be true if we replace homotopy classes over $H/t(G)$ with arbitrary homotopy classes.  As we remarked in Subsection~\ref{subsec:homotopy} above, if $f_0$ and $f_1$ are homotopic maps 
from $M$ into $B$ which are not homotopic over $H/t(G)$ then there is no reason why the structure maps of the induced bundles $f_0^*E$ and $f_1^*E$ should be equal, and hence no reason 
why $f_0^*E$ and $f_1^*E$ should be isomorphic bibundles.  Arbitrary homotopy classes leads to the notion of `concordance' of bibundles: two bibundles $P_0$ and $P_1$ on $M$ are 
said to be {\em concordant} if there is a bibundle $P$ on $M\times I$ such that $P|_{M\times \{0\}}\cong P_0$ and $P|_{M\times \{1\}}\cong P_1$.  The concordance relation is an equivalence 
relation and it can be shown that the set of concordance classes of bibundles on $M$ is in a bijective correspondence with homotopy classes of maps from $M$ into $B$, if $B$ is a space 
as in Theorem~\ref{classif thm}.  So we see that in general the notion of isomorphism of bibundles is a finer equivalence relation, leading to more equivalence classes, than the notion of 
concordance of bibundles.   
\end{note}

As a corollary of Theorem~\ref{classif thm} above, we have the following result.  
\begin{theorem}
\label{thm:explicit class space}
Let $(H,G)$ be a crossed module.  Then the bibundle $EG\times H\to EG\times_G H$ is a classifying bibundle in the sense that there is an isomorphism 
\[ 
\pi_0\bBibun_{(H,G)}(M) \cong [M,EG\times_G H]_{H/t(G)}
\]
for any paracompact space $M$, which is induced by sending the isomorphism class of a bibundle $P$ on $M$ 
to the homotopy class of the map $\Gamma\colon M\to EG\times_G H$ over $H/t(G)$.  
\end{theorem}
\begin{proof}
This follows because $EG\times H\to EG\times_G H$ is a bibundle and $EG\times H$ is contractible as a space over $H$.  
\end{proof}

\subsection{Group structure on the classifying space.}

We will now show that there is a universal $(H, G)$ bibundle $E(H, G)\to B(H, G)$ for which both $E(H, G)$ and $B(H, G)$ are topological groups and the projection map is a group homomorphism.  
Choose a model for the total space $EG$ of the universal $G$ bundle which can be equipped with the structure of a topological group, containing $G$ as a closed subgroup. Assume further  that the action of $H$ on $G$  extends to an action of $H$ on $EG$ by automorphisms. Let us denote this extended action also by 
$$
\a \colon H \to \Aut(EG).
$$
Finally we assume that if $g \in G \subset EG$ then $\a(t(g))(e) = geg^{-1}$. In other words as well
as $G$ being a subgroup of $EG$ we have $(H, G)$ a crossed submodule of $(H, EG)$. 

\begin{example}
Consider the crossed module $\Omega K \to PK$ with the action of $PK$ on $\Omega K$ by conjugation. 
The space $E\Omega K = PK$ is a topological group under pointwise multiplication and the conjugation 
action on $\Omega K$ extends to the adjoint action of $PK$ on  $E\Omega K = PK$.  This is also 
an action by automorphisms. Clearly if $g \in \Omega K$ and $e \in PK$ we have 
$\a(t(g))(e) = \a(g)(e) = g e g^{-1}$. 
\end{example}

More generally, the total space $EG$ of the universal bundle can be constructed as the geometric realization of a certain simplicial space subject to a mild restriction 
on the topological group $G$ (see for instance \cite{May, Ste}).  It turns out that 
that $EG$ carries a natural structure of a topological group, containing $G$ as a closed subgroup.  
It is not hard to show (using the construction of $EG$ given in \cite{May}) that, since $H$ acts by automorphisms on 
$G$, there is an induced action of $H$ on $EG$ and that moreover $\a(t(g))(e) = geg^{-1}$. 

Given this extended action $\a$ we have a  semi-direct product $ EG \rtimes H$ with multiplication
$$
(e, h)(e', h') = (e \a(h)(e') , h h') ).
$$
We denote this group by $E(H, G)$.  Notice that $(e, h)^{-1} = (\a(h)^{-1}(e^{-1}), h^{-1})$.

Notice also that the map  $G \to E(H, G)$  defined by $g \mapsto (g^{-1}, t(g))$ is a homomorphism because 
\begin{align*}
(g^{-1}, t(g) )(k^{-1}, t(k)) 
&= (g^{-1} \a(t(g))(k^{-1}), t(g) t(k))\\
&= (g^{-1} g(k^{-1})g^{-1}, t(g) t(k))\\
 &  = ( k^{-1}g^{-1}, t(gk)) \\
 & = ((gk)^{-1} , t(gk)).
 \end{align*}
Moreover the image  of this map is actually a normal subgroup because
\begin{align}
\label{eq:conjugation}
(e, h) (g^{-1}, t(g)) & (\a(h^{-1})(e^{-1}),  h^{-1}) \\
            & =  (e  \a(h)(g^{-1}), h t(g))             (\a(h^{-1})(e^{-1}), h^{-1}) \notag \\
            & = (e \a(h)(g^{-1}) \a(h t(g))(\a(h^{-1})(e^{-1})), h t(g) h^{-1})  \notag \\
            & = ( (\a(h)(g))^{-1}, t(\a(h)(g)) ). \notag 
            \end{align}
                       
We denote the image of this homomorphism inside $E(H, G)$ by $G$ and hence have the exact 
sequence of topological groups
\begin{equation}
\label{eq:g>B(H, G)>E(H, G)}
G \to E(H, G) \to B(H, G).
\end{equation}
In fact the map $E(H, G)\to B(H, G)$ admits local sections since it is obtained by pullback from the universal $G$ bundle $EG\to BG$ via the projection $B(H, G)\to BG$ (see \cite{BaeSte} for a proof of this).  
Notice that right action by $G$ on $E(H,G)$  is 
$$
(e, h) g = (e, h)(g^{-1}, t(g) ) = (e\a(h)(g^{-1}), h t(g)).
$$
Example~\ref{ex:groups} shows that $E(H,G)\to B(H,G)$ is a $(G,E(H,G))$ bibundle, where the structure map 
$E(H,G)\to E(H,G)$ is the identity.  Notice though that there is a morphism of crossed modules 
$(G,E(H,G))\to (G,H)$, the homomorphism $E(H,G)\to H$ being the projection onto $H$ in the semi-direct product.   
Therefore there is a natural extension of $E(H,G)\to B(H,G)$ to an
 $(H,G)$ bibundle  (see the discussion in subsection~\ref{sec:ext-red}).  
The structure map $\Phi\colon E(H,G)\to H$ is then given by $\Phi(e,h) = h$.  

Consider now the bijection $ \hat \chi \colon EG \times H \to E(H, G)$ defined by
$$
\begin{array}{ccc}
\hat \chi \colon EG \times H & \to & E(H, G) \\
(e, h) & \mapsto & (\a(h)(e^{-1}), h)
\end{array}
$$
and note that $\chi$ commutes with structure maps which are the projections
onto $H$ in both cases.  To see that $\chi$ is a bibundle
isomorphism we only have to check that it commutes with the right action of $g \in G$.
We have
\begin{align*}
\hat\chi((e, h) g) &= \hat\chi( (eg, h t(g)) ) \\
     & = ( \a(h t(g))( g^{-1} e^{-1}), h t(g) ) \\
      & = ( \a(h)(\a( t(g))( g^{-1} e^{-1})), h t(g) ) \\
       & = ( \a(h)(g (g^{-1} e^{-1})g^{-1}), h t(g) ) \\
       & = ( \a(h)( e^{-1}g^{-1}), h t(g) ) \\
           & = ( \a(h)( e^{-1}) \a(h)(g^{-1}), h t(g) ) \\
           &= (\a(h)(e^{-1}), h)(g^{-1}, t(g) )\\
            & = \hat\chi(e,h) g.
     \end{align*}
      
     It follows that $\hat\chi$ induces a map  $\chi \colon EG \times_G H \to B(H, G) $
     and that

\begin{equation} 
\label{} 
\xy 
(-15,7.5)*+{EG \times H}="1"; 
(10,7.5)*+{E(H, G)}="2"; 
(-15,-7.5)*+{EG \times_G H}="3"; 
(10,-7.5)*+{B(H, G)}="4"; 
{\ar^-{\hat\chi} "1";"2"}; 
{\ar "1";"3"}; 
{\ar_-{\chi} "3";"4"}; 
{\ar "2";"4"};
\endxy
\end{equation}
is an isomorphism of bibundles.

The diagram~\eqref{universal commuting diagram} above now becomes a diagram of topological groups and continuous homomorphisms between them 
\begin{equation*}  
\xy 
(-10,7.5)*+{E(H, G)}="1"; 
(10,7.5)*+{H}="2"; 
(-10,-7.5)*+{B(H, G)}="3"; 
(10,-7.5)*+{H/t(G).}="4"; 
{\ar^-{\Psi} "1";"2"}; 
{\ar "1";"3"}; 
{\ar_-{\Phi} "3";"4"}; 
{\ar "2";"4"};
\endxy
\end{equation*}
If $P \to M$ is a bibundle we define $F \colon M \to B(H, G)$ by composing 
 the function $\Gamma \colon M \to EG \times_G H$ with the function $\chi$ to obtain 
 the composite diagram
 
 \begin{equation*} 
\xy 
(-15,7.5)*+{P}="1"; 
(10,7.5)*+{EG \times H}="2"; 
(35,7.5)*+{E(H, G)}="3"; 
(-15,-7.5)*+{M}="4"; 
(10,-7.5)*+{EG \times_G H}="5"; 
(35,-7.5)*+{B(H, G)}="6"; 
{\ar^-{\hat \Gamma} "1";"2"}; 
{\ar^-{\hat \chi} "2";"3"}; 
{\ar_-{\Gamma} "4";"5"}; 
{\ar_-{\chi} "5";"6"};
{\ar "1";"4"}; 
{\ar "2";"5"}; 
{\ar "3";"6"}; 
\endxy
\end{equation*}
from which we obtain 

 \begin{equation*}  
\xy 
(-15,7.5)*+{P}="1"; 
(10,7.5)*+{E(H, G)}="2"; 
(35,7.5)*+{H}="3"; 
(-15,-7.5)*+{M}="4"; 
(10,-7.5)*+{B(H, G)}="5"; 
(35,-7.5)*+{H/t(G).}="6"; 
{\ar^-{\hat F} "1";"2"}; 
{\ar^-{\Psi} "2";"3"}; 
{\ar_-{F } "4";"5"}; 
{\ar_-{\Phi} "5";"6"};
{\ar "1";"4"}; 
{\ar "2";"5"}; 
{\ar "3";"6"}; 
\endxy
\end{equation*}
Clearly $E(H,G)\to B(H,G)$ is a classifying bibundle for $(H,G)$ bibundles.  Since $E(H,G)\to B(H,G)$ is a bibundle arising 
from the quotient of groups~\eqref{eq:g>B(H, G)>E(H, G)}, 
it has a nice behaviour with respect to products of bibundles.  More precisely we have the following results.   
The discussion in Example~\ref{ex: group bibundle morphism} gives us the following lemma.
\begin{lemma} 
\label{products of bibundles} \hfil
\begin{enumerate}
\item   The product in the group $E(H, G)$ induces a morphism of bibundles $E(H, G)\otimes E(H, G)\to E(H, G)$ covering 
the product $B(H, G)\times B(H, G)\to B(H, G)$ in the group $B(H, G)$.  
\item The inverse map $B(H, G) \to B(H, G)$ pulls back the bibundle $E(H, G) \to B(H, G)$ to its dual.
\end{enumerate}
\end{lemma}
From this we easily deduce the next proposition.
  \begin{proposition}
  \label{mult prop of class map}\hfil
 \begin{enumerate}
 \item 
  Let $F_1, F_2 \colon M \to B(H, G)$ and define $F_1F_2 \colon M \to B(H, G)$ to be the pointwise product. Then 
  $(F_1F_2)^*(E(H, G)) \simeq F_1^*(E(H, G)) \otimes F_2^*(E(H, G))$. 
  \item Let $F \colon M \to B(H, G)$ and denote by $F^{-1}$ the pointwise inverse. Then $(F^{-1})^*(E(H, G)) \simeq (F^*(E(H, G)))^* $.
  \end{enumerate}
    \end{proposition}
As we have already observed, $E(H,G)\to B(H,G)$ is a universal bibundle in the sense of Theorem~\ref{classif thm}
and hence there is an isomorphism 
\begin{equation}
\label{univ isom for prod}
[M,B(H,G)]_{H/t(G)}\cong \pi_0\bBibun_{(H,G)}(M).
\end{equation}
Since $B(H,G)$ is a topological group the set $[M,B(H,G)]_{H/t(G)}$ of homotopy classes of maps 
over $H/t(G)$ acquires a natural structure as a topological group.  As we have remarked previously 
(see Lemma~\ref{group structure on IBibun}) 
 $\pi_0\bBibun_{(H,G)}(M)$ also has a natural structure of a group where the product $[P]\cdot [Q]$ 
is the isomorphism class of the bibundle $P\otimes Q$.  Proposition~\ref{mult prop of class map} shows 
that the isomorphism~\eqref{univ isom for prod} preserves products, since 
\[
[f^*E(H,G)]\cdot [g^*E(H,G)] = [(fg)^*E(H,G)].  
\]
It follows that~\eqref{univ isom for prod} is an isomorphism of groups.  We record this observation in the following theorem.  
\begin{theorem} 
Let $(H,G)$ be a crossed module.  Then the bibundle $E(H,G)\to B(H,G)$ is a classifying bibundle which preserves 
group structures in the sense that there is an isomorphism of groups 
\[
\pi_0\bBibun_{(H,G)}(M)\cong  [M,B(H,G)]_{H/t(G)}
\]
for any paracompact space $M$.  
\end{theorem}   
    
\subsection{$(H, G)$ bibundle structures on a $G$ principal bundle.}
Consider the commuting diagram
\begin{equation*} 
\xy 
(10,7.5)*+{EG \times H}="2"; 
(35,7.5)*+{EG}="3"; 
(10,-7.5)*+{EG \times_G H}="5"; 
(35,-7.5)*+{BG}="6"; 
{\ar_-{\pi_{BG}} "5";"6"};
{\ar_-{\pi_{EG}} "2";"3"};
{\ar "2";"5"}; 
{\ar "3";"6"}; 
\endxy
\end{equation*}
where $\pi_{EG}$ and $\pi_{BG}$ are the natural projections. This is a morphism of $G$ bundles
and  shows  that $\pi_{BG}$ is a classifying map for the $G$ bundle
$EG \times H \to EB \times_G H$.  
Using the isomorphism $\chi$ from the bibundle $EG \times H \to EG \times_G H$  to the bibundle $E(H, G) \to B(H, G)$, 
which is also a $G$ bundle isomorphism, 
we see that $\pi \colon B(H, G) \to BG$  defined by $\pi = \pi_{BG} \circ \chi^{-1} $ is a 
 classifying map for the $G$ bundle $E(H, G) \to B(H, G)$. Note that $\pi \colon B(H, G) \to BG $, 
 like $\pi_{BG}$, has fibre $H$. 

It follows that if  $P \to M$ is an $(H, G)$ bibundle with classifying map $F \colon M \to B(H, G)$ then $\pi \circ F \colon M \to BG$ is a classifying map for the $G$ bundle $P \to M$. Conversely if  a $G$ bundle   $P \to M$ has a classifying map  $f \colon M \to BG$ which lifts to a map $\hat f  \colon M \to B(G, H)$ then it is isomorphic to $f^*(E(H, G))$
and thus admits an $(H, G)$  bibundle structure.  Hence a $G$ bundle  $P \to M$ admits an ($H, G)$
bibundle structure if and only if it has  a classifying map $M \to BG $ which lifts to $M \to B(H, G)$. 
  Thus we have the following proposition. 
    
\begin{proposition} 
    If $H$ is contractible then every $G$ bundle admits an $(H, G)$ bibundle structure.
    \end{proposition}
    
    \begin{example} Consider a $\Omega K$ bundle $P \to M$.  Then as $PK$ 
    is contractible we can always lift a map $M \to B\Omega K $ to $B(PK, \Omega K)$ and every
    $\Omega K$ admits a $(PK, \Omega K)$ bibundle structure. Recall that a $(PK, \Omega K)$ is also
    an $\Omega K$ bibundle so every $\Omega K$ bundle admits an $\Omega K$ bibundle structure.
    \end{example}

\subsection{Loop groups} 
We have seen that the existence of  bibundles that are not abelian relates to the size of $\Out(G)$ 
and we have commented that we are therefore interested in groups that have 
large outer automorphism groups. One example is the group 
$G = \Omega K$ of based loops in a compact Lie group $K$.  
There are a number of possible groups $H$ and homomorphisms making $\Omega K \to H$
into a crossed module so for the moment we will make some general comments, before looking at some specific examples.
 
As we remarked above we have $E\Omega K = PK$ and $B\Omega K = K$ with the projection 
$PK \to K$ being evaluation at $1$. 
We have that $EG \times H = PK \times H$ where the action of $
\Omega K$ is $(e, h)k = (ek, h t(k))$ and that $E(H, \Omega K) = PK \times H $ with the 
product $(e, h)(e', h') = (e\a(h)(e'), h h')$ and the action of 
$k \in \Omega K$ being $(e, h)k = (e\a(h)(k^{-1}), h t(k))$. The projection 
from $B(H, \Omega K)$ to $B\Omega K = K$ is given by 
$[e, h] \mapsto \a(h)^{-1}(e^{-1})(1)$. The type map $B(H, \Omega K) \to H/t(\Omega K)$ is just
projection of $(e, h)$ to the coset of $h$ in $H/t(\Omega K)$.

\begin{example} 
\label{loop group example}
Let $H= PK$ with the action $\a(h)(e) = h e h^{-1}$ using pointwise multiplication of paths. 
Then  $E(PK, \Omega K) = PK \times PK$ with the multiplication 
$$
(e, h)(\bar e, \bar h) = (e h \bar e h^{-1}, h \bar h).
$$
We can identify this with the usual product $PK \times PK$ by the 
isomorphism $(e, h ) \mapsto (eh, h)$. The subgroup $\Omega K$ of all $(k^{-1}, k)$
becomes the subgroup of all $(1, k)$ or $\{ 1 \} \times \Omega K$. So we have 
$E(PK, \Omega K) = PK \times PK$ and $B(PK, \Omega K) = PK \times K$.  The projection of 
a pair $(e, k) \in PK \times K $ to $B\Omega K = K$ can be calculated by first 
reversing the isomomorphism above to send it to $(ek^{-1}, k)$ and then mapping this to $\a(k)^{-1}((e k^{-1})^{-1})(1) = e(1)k$. 
Notice that any map $ f \colon M \to K$ can be lifted to $F \colon M \to PK \times K$ by taking $F(m) = (1, f(m))$. 

In this example $H / t(\Omega K) = PK/\Omega K = K$ and the type map $B(PK, \Omega K) = PK \times K  \to K$ is just projection onto the second 
factor. 
It follows that we have the exact sequence of groups
$$ 
1\to PK\to B(PK, \Omega K) \to K \to 1.  
$$
This is an instance of an exact sequence of groups 
\[
1\to EG/G_1\to B(H,G)\to H/t(G)\to 1 
\]
which exists for any crossed module $(H,G)$.  
\end{example}    

Using these observations we illustrate how the set of homotopy classes of maps $[M,B(H,G)]$ can fail to be isomorphic 
to the set of isomorphism classes of $(H,G)$-bibundles, as mentioned earlier.  

\begin{example} 
\label{counterexample}
Let $K$ be a compact connected, non-trivial Lie group and 
let $G = \Omega K$.  Then we have seen in Example~\ref{loop group example} above 
that $B(\Omega K, PK) $ is isomorphic to the group 
$PK \times K$ and the homomorphism $\Phi\colon B(\Omega K , PK) \to K$ is the natural projection.  
Clearly we can find homotopic maps $F_0,F_1\colon M\to PK \times K$ 
such that $\Phi F_0$ and $\Phi F_1$ are not equal.  Hence the 
$\Omega K$ bibundles induced by pull back with $F_0$ and $F_1$ are not isomorphic.  
\end{example}

We leave it for the interested reader to consider some of the other crossed modules $(H, \Omega K)$ in the following 
examples. 

\begin{example} The group  $\Diff_0([0,1])$ of diffeomorphisms fixing $0$ and $1$ acts on $\Omega K$ so we can take $H$ equal to the 
semi-direct product $PK \rtimes \Diff_0(S^1)$.
\end{example}

\begin{example} 
Replace $PK$ by the group of  smooth maps from $[0, 1]$ to $\Aut(K)$ acting pointwise. As automorphisms
fix the identity there is no need to impose a condition on the map at the endpoints. 
\end{example}

\section{Conclusion}  While bibundles are of independent interest one motivation for discussing them is the notion 
of an $(H, G)$ 
bibundle gerbe \cite{AshCanJur}. We will indicate here briefly how our approach 
will apply in this case. A complete discussion will appear in \cite{MurRobSte}. 

We assume the reader is familiar with abelian bundle gerbes \cite{Murray}.  First we have the analogue of the definition of an abelian bundle gerbe.

\begin{definition}[c.f. \cite{AshCanJur}]
An {\em $(H, G)$ bibundle gerbe} on $M$, or just bibundle gerbe when the crossed module $(H,G)$ is understood, consists of a 
pair $(P,Y)$ where $\pi\colon Y\to M$ is a surjective submersion and 
$P\to Y^{[2]}$ is an $(H, G)$ bibundle equipped with a bibundle 
gerbe product.   This is a bibundle map which on fibres takes the form 
$$ 
P_{(y_1,y_2)}\otimes P_{(y_2,y_3)}\to P_{(y_1,y_3)}
$$ 
for $(y_1,y_2,y_3)\in Y^{[3]}$.  The bibundle gerbe product is required to be associative in the usual sense.  
\end{definition}

Fundamental to the theory of bibundle gerbes is the notion of stable isomorphism.  Before we give the definition, observe that if $(P,Y)$ is a bibundle gerbe and $R$ is a bibundle on $Y$, then 
we can construct a new  bibundle gerbe $(\pi_2^*R^*\otimes P\otimes \pi_1^*R,Y)$ on $M$ with bibundle 
gerbe product given fibrewise by using the bibundle gerbe product on $P$ and 
contraction as 
$$ 
R_{y_1}^*\otimes P_{(y_1,y_2)}\otimes R_{y_2}\otimes R_{y_2}^*\otimes P_{(y_2,y_3)}\otimes R_{y_3}\to R_{y_1}^*\otimes P_{(y_1,y_3)}\otimes R_{y_3}.  
$$ 
With this construction in hand we can make the following definition, which is closely related to Definition 12 in \cite{AshCanJur}.  
\begin{definition} 
\label{def:stable-bibundle-gerbes}
Let $(P,Y)$ and $(P',Y')$ be bibundle gerbes on $M$.  We say that $P$ is {\em stably isomorphic} to $P'$ if 
there exists a bibundle $R$ on $Y\times_M Y'$ together with an isomorphism 
of bibundle gerbes 
\begin{equation} 
\label{eq: stable iso} 
\pi_2^*R^*\otimes P\otimes \pi_1^*R\cong P' 
\end{equation}
where we have suppressed the projections $(Y\times_M Y')^{[2]}\to Y^{[2]}$ and $(Y\times_M Y')^{[2]}\to (Y')^{[2]}$ which are used to pullback $P$ and $P'$ respectively.  
\end{definition} 

Finally the type of an $(H, G)$  bibundle gerbe is defined as follows.
Let $(P, Y)$ be an $(H, G)$  bibundle gerbe over $M$. Then $P \to Y^{[2]} $ has a type
map $Y^{[2]} \to H/t(G)$. The existence of the bundle gerbe multiplication and its
associativity means that 
\begin{equation}
\label{prebundle cocycle equation}
\phi(y_1,y_2)\phi(y_2,y_3) = \phi(y_1,y_3)
\end{equation}
for all $(y_1, y_2, y_3) \in Y^{[3]}$. 
Recall \cite{CarJohMurWan} that a $K$ {\em prebundle} over $M$ consists of a submersion $Y \to M$ and a map 
$k \colon Y^{[2]} \to K$  satisfying a cocycle equation analogous to~\eqref{prebundle cocycle equation} above.   
Every $K$ prebundle over $M$ determines a principal $K$-bundle over $M$ and conversely.  

It follows that every $(H,G)$ bibundle gerbe $(P, Y)$ over $M$ defines an $H/t(G)$ prebundle
$(\phi, Y)$ over $M$ and hence a principal $H/t(G)$ bundle over $M$.  
We call this  the {\em type} of the bibundle gerbe.  In the sheaf theoretic setting, this pre-bundle is 
known as the {\em band} of the gerbe.

\begin{example}
A Jandl bundle gerbe \cite{NikSch} consists of a $(\ZZ_2, U(1))$ bundle gerbe $(P, Y)$
over $M$.  The type of $(P, Y)$ is a $\ZZ_2$ prebundle over $M$ and the induced $\ZZ_2$-bundle
is called the {\em orientation bundle} of the Jandl bundle gerbe.
\end{example}
We will discuss all of these notions in more detail in \cite{MurRobSte}.

\end{document}